\documentclass{article}

\usepackage{amssymb}
\usepackage{dsfont}

\newtheorem{lem}{Lemma}[section]
\newtheorem{defi}[lem]{Definition}
\newtheorem{theo}[lem]{Theorem}
\newtheorem{prop}[lem]{Proposition}

\newtheorem{rem}[lem]{Remark}

\newtheorem{fait}[lem]{Fact}
\newtheorem{cor}[lem]{Corollary}
\newtheorem{claim}[lem]{Claim}
\newtheorem{question}[lem]{Question}

\newtheorem{example}[lem]{Example}

\newcommand{\proof}{\noindent{\bf Proof.}~}
\newcommand{\qed}{\ \hfill$\square$\bigskip}

\newcommand{\ad}{\hbox{\rm ad}\,}
\newcommand{\eq}{\hbox{\rm eq}\,}

\newcommand{\PSL}{\hbox{\rm PSL}\,}

\newcommand{\Tor}{\hbox{\rm Tor}\,}

\newcommand{\<}{\langle}

\renewcommand{\>}{\rangle}
\renewcommand{\o}{^\circ}

\newcommand{\Q}{\mathds{Q}}

\newcommand{\N}{\mathds{N}}
\newcommand{\R}{\mathds{R}}

\title{Commutators in groups definable in o-minimal structures}
\author{El\'ias Baro, Eric Jaligot, and Margarita Otero\thanks{Corresponding author:  e-mail (margarita.otero@uam.es); tel.(34 914973808); FAX(34 914974889). The first and third authors are partially supported by MTM2008-00272 and Grupos UCM 910444. } 
\\ Departamento de Matem\'aticas \\ Universidad Aut\'onoma de Madrid\\ 28049 Madrid, Spain}
\date{May 11, 2010}

\begin{document}
\maketitle
\begin{abstract}We prove the definability, and actually the finiteness of the commutator width, of many 
commutator subgroups in groups definable in o-minimal structures. It applies in particular to derived series 
and to lower central series of solvable groups. Along the way, we prove some 
generalities on groups with the descending chain condition on definable subgroups 
and/or with a definable and additive dimension. 
\end{abstract}

{\em MCS2010}: 03C64; 20F12, 20F38, 20A15, 03C60.

\section{Introduction}

The development of model theoretic stability since the 70's provided an effective bridge from 
classical geometric objects to general first-order definable sets. 
For example, the voluminous theory of groups of finite Morley rank can be seen as a 
large generalization of that of algebraic groups over algebraically closed fields. More recently, 
the theory of groups definable in o-minimal structures has had important developments as well, 
often recovering some aspects of stable groups theory. The present paper continues in this vein. 

Groups definable in o-minimal structures and groups of finite Morley rank share many properties. 
They are both equipped with a finite dimension which is definable and additive, satisfy 
the descending chain condition on definable subgroups, and have definably connected 
components, i.e., smallest definable subgroups of finite index. 
These properties suffice for many developments of the theory of 
groups of finite Morley rank. However, if one restricts our attention to the 
natural examples of groups definable over fields of characteristic $0$, 
the category of groups definable in o-minimal structures appears to be strictly larger since it contains 
all semi-algebraic groups over real closed fields in addition to algebraic groups over 
algebraically closed fields of characteristic $0$. 

In fact, the main difference between groups of finite Morley rank and groups 
definable in an o-minimal structure is in the behaviour of the dimension. 
The dimension of a definable set $A$ in the o-minimal case fails the main property 
of the Morley rank: 
$\dim(A)\geq n+1$ if and only if $A$ contains infinitely many pairwise disjoint 
definable subsets $A_i$ with $\dim(A_i)\geq n$. 
This essential definition of the Morley rank is crucial in Zilber's stabilizer argument, and consequently 
Zilber's generation lemma on indecomposable sets in the finite Morley rank context \cite[Section 5.4]{MR1321141}. And this is what allows one to prove the definability of many commutators subgroups, 
and in particular of derived subgroups, in this context. 

For groups definable in o-minimal structures, derived subgroups need not be definable. 
Using recent results from \cite{HPPCentExt} on central extensions in the o-minimal case, 
Annalisa Conversano exhibits in \cite[Example 3.1.7]{ConversanoThesis} 
a definably connected group $G$ definable in an o-minimal expansion of the reals 
with $G'$ not definable. Furthermore, $G$ is a central extension, by an infinite 
center isomorphic to $[0,1[$, of the simple group $\PSL_2(\R)$. In the present paper we will essentially 
prove that central extensions of this type are the only obstructions to the definability 
of commutator subgroups. In order to state our main theorem we 
consequently need the following definition. 

\begin{defi}\label{DefiStrictCentralExt}
We say that a definably connected group $G$ definable in an o-minimal structure 
is a {\em strict} central extension of a definably simple group if $Z(G)$ is infinite and 
$G/Z(G)$ is infinite nonabelian and definably simple. 
\end{defi}

Definability, in Definition \ref{DefiStrictCentralExt} and throughout the paper, always refers 
to a fixed structure, typically an o-minimal structure. 

We call {\em section} of a group $G$ any quotient $H/K$ where $K \trianglelefteq H \leq G$, 
and we speak of {\em definable} section when both $K$ and $H$ are definable. 
In view of the previous example of a group definable in 
an o-minimal structure and with a non definable derived subgroup, we will mostly work with the following 
assumption. It consists merely in excluding ``bad" sections of this type. 

\begin{defi}\label{assumption}
We say that a group $G$ definable in an o-minimal structure satisfies the {\em assumption (*)} whenever 
the derived subgroup $(H/K)'$ is definable for every definable section $H/K$ of $G$ which is a strict central extension of a definably simple group. 
\end{defi}

Our main theorem is the following. 

\begin{theo}\label{maintheo}
Let $G$ be a group definable in an o-minimal structure, and 
$A$ and $B$ two definable subgroups which normalize each other  and 
such that $A\o B\o$ satisfies the assumption $(*)$. 
Then the subgroup $[A,B]$ is definable and $[A,B]\o=[A\o,B][A,B\o]$. 
Furthermore, any element of $[A,B]\o$ can be expressed as the product of at most 
$\dim([A,B]\o)$ commutators from $[A\o,B]$ or $[A,B\o]$ whenever $A\o$ or $B\o$ is solvable. 
\end{theo}

We will see in Corollary \ref{CorABSolvAB*} below that if $A\o$ and $B\o$ are solvable in 
Theorem \ref{maintheo}, then $A\o B\o$ satisfies the assumption $(*)$ and thus $[A,B]$ is always definable 
in this case. In particular, Theorem \ref{maintheo} implies that if $G$ is a group definable in an 
o-minimal structure with $G\o$ solvable, then the lower central series and derived series of $G$ consist 
of definable subgroups (see Corollary \ref{CorGoSolvGnGndef} below). 
In view of Conversano's example on the other hand, the assumption $(*)$ is necessary if one wants a 
statement for definability as general as in Theorem \ref{maintheo}, 
though of course the conclusion of the theorem may remain valid without that assumption 
in more specific cases. 

A commutator group $[A,B]$ as in Theorem \ref{maintheo} is the countable union of the definable sets 
$[A,B]_n$, consisting of the products of at most $n$ commutators $[a,b]$, or their inverses $[a,b]^{-1}$, 
with $a$ in $A$ and $b$ in $B$. Our proof of Theorem \ref{maintheo} will consist in showing that $[A,B]$ has 
{\em finite commutator width}, i.e., that $[A,B]=[A,B]_n$ for some $n$. This is equivalent to the 
definability of $[A,B]$ when the ground structure in which $G$ is defined 
is $\omega$-saturated, but since we make no saturation assumption 
we show directly the finiteness of the commutator width. Besides, we will always try to keep 
the best control on the commutator width, generally in terms of the dimension of the groups involved. 
We emphasize that even the existence of such a finite $n$ did not seem to be known, even in 
the more classical case of semi-algebraic groups over real closed fields. Our result, when it applies, 
also says something apparently new about the commutator width in Lie groups. 

Our argument for the proof of Theorem \ref{maintheo} will consist mainly in finding very 
rudimentary forms of Zilber's general stabilizer argument on generation by indecomposable sets 
in the finite Morley rank context. In fact, our bounds $n$ as above we will be obtained ultimately 
with Corollary \ref{Lem[X,H]definable} below, concerning the generation by definable and 
definably connected subgroups in abelian groups. 

We refer to \cite{MR1633348} for a general introduction to o-minimal structures, 
and to \cite{MR2436142} for groups definable in o-minimal structures. We will 
frequently point out analogy to the finite Morley rank case in 
\cite{MR1321141}. All model theoretic notions used here are rather elementary and can be found in 
any introductory book on model theory. We simply recall that an {\em o-minimal} structure 
is a first-order structure $\cal M$ with a total, dense, and without end-points definable order 
and such that every definable subset of $M$ is a boolean combination of intervals 
with end-points in $M\cup\{\pm\infty\}$.  

The paper is organized as follows. In Section \ref{SectionCommutatorsPrelude} we recall or 
develop the background used on commutators in arbitrary groups. 
In Sections \ref{Gpsdcc} and \ref{NilpGpsdcc} we take the opportunity to recast, when this is possible, 
various facts typical of groups of finite Morley rank to the more general context of groups 
with various descending chain conditions on definable subgroups. 
In Section \ref{SectionDefAddDimension} 
we continue with an axiomatic treatment of groups with a definable and additive dimension. 
In Section \ref{Section6} 
we insert groups definable in o-minimal structures in the preceding abstract machinery. 
Section \ref{Section7} 
is devoted to a preliminary analysis of the structure of solvable groups. 
In Section \ref{Section8} 
we prove our main results about commutators and finally in 
Section \ref{Section9} 
we conclude with further results and questions. 

\medskip
We thank Alessandro Berarducci, Gregory Cherlin, Ya'acov Peterzil, and Anand Pillay for 
various related conversations at Oberwolfach in January 2010. We also thank 
Dragomir Djokovic for hints concerning commutators in simple real Lie groups. 

\section{Prelude to commutators}\label{SectionCommutatorsPrelude}

We fix the following notation as far as commutators are concerned in groups. 
If $a$ and $b$ are two elements of a group $G$, then we let $a^b=b^{-1}ab$ and $[a,b]=a^{-1}a^b$. 
If $A$ and $B$ are two arbitrary subsets of $G$, then 
$[A,B]_n$ denotes, for $n\geq 1$, the set of products of at most $n$ elements of the form 
$[a,b]$ or $[a,b]^{-1}$, with $(a,b)$ in $A\times B$. We denote by $[A,B]_{\infty}$, or simply by 
$[A,B]$ when there is no risk of confusion, the group 
generated by all commutators $[a,b]$, that is the union of all sets $[A,B]_n$. The 
{\em commutator width}, also frequently called the {\em commutator length}, of $[A,B]$ is the smallest $n$ 
such that $[A,B]=[A,B]_n$ if such an $n$ exists, and infinite otherwise. 

Since our proofs by induction will use passages to quotients we will frequently ``lift" 
commutators. This simply means that if $N$ is a normal subgroup of $G$ and $(a,b)\in G^2$, then 
$[aN,bN]=[a,b]N$ in $G/N$. In order to keep track of the commutator width when passing to quotients, 
the following lemma will be helpful. 

\begin{lem}\label{LemLiftCommutators}
Let $G$ be a group, $N$ a normal subgroup of $G$, and $A$ and $B$ two subsets of $G$ such 
that, modulo $N$, $[A,B]=[A,B]_k$ for some $k$. 
Suppose also that $[A,B]\cap N=[A,B]_s\cap N$ for some $s$. 
Then $[A,B]=[A,B]_{k+s}$. 
\end{lem}
\proof
Clearly, $[A,B]N/N=[AN/N,BN/N]$ by lifting of commutators. Our assumption on the commutator width of 
$[A,B]$ modulo $N$ reads as 
$$[AN/N,BN/N]=[AN/N,BN/N]_k.$$

Take any element $g$ in $[A,B]$. Since $g$ is also in $[A,B]N$, $g=\alpha u$ for some 
$\alpha$ in $[A,B]_k$ and some $u$ in $N$. But now $u=\alpha^{-1}g$ is in $[A,B]\cap N$, and by assumption 
$u$ is in $[A,B]_s$. It follows that $g=\alpha u$ is in $[A,B]_{k+s}$. 
\qed

\bigskip
The following consequence of a theorem of Baer will be used for dealing with non-connected 
``finite bits" in our proofs. 

\begin{lem}\label{CorBaer}
Let $A$ and $B$ be two subgroups of a group, with $A$ normalizing $B$. Suppose 
$C_A(B)$ of finite index in $A$ and $C_B(A)$ of finite index in $B$. 
Then $[A,B]$ is a finite subgroup of $B$. 
\end{lem}
\proof
Let $n$ be the index of $C_A(B)$ in $A$ and $m$ the index of $C_B(A)$ in $B$. Fix 
$a_i$, $1\leq i \leq n$, a set of representatives in $A$ of $A/C_A(B)$, and 
$b_j$, $1\leq j \leq m$, a set of representatives in $B$ of $B/C_B(A)$. 

We claim that the set of commutators $\{[a,b]~|~(a,b)\in{A\times B}\}$ is finite. 
In fact, any element 
$a$ in $A$ has the form $a_i\alpha$ for some $i$ and some $\alpha$ in $C_A(B)$ and 
any element $b$ in $B$ has the form $b_j\beta$ for some $j$ and some $\beta$ in $C_B(A)$. 
Now $[a,b]=[a_i\alpha,b_j\beta]={\alpha}^{-1}{a_i}^{-1}{\beta}^{-1}{b_j}^{-1}a_i{\alpha}b_j\beta$. 
Since $\alpha$ centralizes $B$ and $\beta$ centralizes $A$, we get 
$[a,b]={\beta}^{-1}{\alpha}^{-1}[a_i,b_j]{\alpha}{\beta}=[a_i,b_j]^{{\alpha}{\beta}}$, and since 
$[A,B]\leq {B}$ we also get 
$$[a,b]=[a_i,b_j]^{\beta}=[a_i,b_j^{\beta}].$$
Hence $[a,b]\in a_i^{-1}{a_i}^B$. But $C_B(A)\leq C_B(a_i)$ and it follows that $a_i$ has a centralizer 
of finite index in $B$, and in other words that ${a_i}^B$ is finite. Hence $[a,b]$ varies in a finite set, and this 
proves our claim. 

Now it suffices to apply a result of Baer saying that if $A$ and $B$ are two subgroups of a group, 
with $A$ normalizing $B$ and such that the set of commutators $\{[a,b]~|~(a,b)\in{A\times B}\}$ is finite, 
then the group $[A,B]$ is finite \cite{MR0147553}. 
\qed

\bigskip
We recall two classical formulas which are valid for any elements $a$, $b$, and $c$ in a group $G$: 
$$[a,bc]=[a,c][a,b]^c~\mbox{and}~[ab,c]=[a,c]^b[b,c].$$
With these two formulas one sees that if $A$ and $B$ are two arbitrary subgroups of a group $G$, then 
the subgroup $[A,B]$ is always normalized by $A$ and $B$. The proof of our main theorems 
will reduce in a crucial way to the second of these formulas via the following lemma. 

\begin{fait}\label{FactGenAdMap}
Let $G$ be a group, $x$ an element of $G$, and $H$ a subgroup such that 
${\{[h,x]~|~h\in H\}}\subseteq C_G(H)$. Then the map 
$$\begin{array}{llll}
    \ad_x : & H & \longrightarrow & G\\
                 & h  & \longmapsto      & [h,x].
\end{array}$$
is a group homomorphism, with $C_H(x)$ as kernel and $\{[h,x]~|~h\in H\}$ as image. 
\end{fait}
\proof
For any $h_1$ and $h_2$ in $H$ we have 
$$[h_1h_2,x]=[h_1,x]^{h_2}[h_2,x]=[h_1,x][h_2,x]$$
by assumption. This shows that $\ad_x$ is a group homomorphism and the rest is clear. 
\qed

\bigskip
For any group $G$ we let $G^1=G^{(1)}=G'=[G,G]$, the {\em derived} subgroup of $G$, and we define 
by induction on $n\geq 1$ the {\em lower central series} and {\em derived series} as follows: 
$G^{n+1}=[G,G^n]$ and $G^{(n+1)}=[G^{(n)},G^{(n)}]$. We say that $G$ is {\em $n$-nilpotent} 
(resp. {\em $n$-solvable}) whenever $G^n=1$ (resp. $G^{(n)}=1$), and we say that $G$ is {\em nilpotent} 
(resp. {\em solvable}) whenever it is $n$-nilpotent (resp. $n$-solvable) for some $n\geq 1$. The nilpotency 
(resp. solvability) {\em class} of $G$ is the smallest $n$ such that $G$ is $n$-nilpotent (resp. $n$-solvable). 
We also define the {\em upper central series} as follows: $Z_1(G)=Z(G)$ is the center of $G$, and by induction 
$Z_{n+1}(G)$ is the preimage in $G$ of $Z(G/Z_n(G))$. As far as nilpotent groups are concerned, our 
inductive proofs will naturally use the following elementary fact. 

\begin{fait}\label{FactZnGnNilpGps}Let $n\geq 1$. 
\begin{itemize}
\item[$(a)$]
If $G$ is a nilpotent group of class $n$, then $G^i\leq Z_{n-i}(G)$, for $0\leq i \leq n$. 
\item[$(b)$]
A group $G$ is nilpotent of class $n$ if and only if $Z_n(G)=G$ and $Z_{n-1}(G)<G$. In this case $G/Z(G)$ 
is nilpotent of class $n-1$. 
\end{itemize}
\end{fait}
\proof
The inclusions of first claim are easily proved by induction on $n$; see for 
instance \cite[Lemma 0.1.7]{MR1473226}. 
See \cite[Corollary 0.1.8]{MR1473226} for the first point of the second claim; the second point then follows. 
\qed

\section{Groups with $dcc$}\label{Gpsdcc}

In this section we work with a group definable in a structure $\cal M$ under the 
mere assumption that it satisfies the {\em descending chain condition}, $dcc$ for short, 
on definable subgroups: any strictly descending chain of $\cal M$-definable subgroups is finite.
Groups definable in an o-minimal structure $\cal M$ satisfy the $dcc$ on definable subgroups 
\cite[Remark 2.13 (ii)]{MR961362}, 
and consequently all the results of the present section will apply to such groups. 

In what follows definability always refers to definability from the original universe $\cal M$. 
Actually everything in the present section can be done for a group interpretable in $\cal M$, and hence 
definability here may refer to definability in ${\cal M}^{\eq}$. 
Any definable group with no proper definable subgroup of finite index is said to be {\em definably connected}.  
Notice that a definable group can be definably connected even without satisfying the $dcc$. We start 
with two general facts about definably connected groups (independent of the $dcc$). 

\begin{fait}\label{ActionConOnFinite}
Let $G$ be a definably connected group. 
\begin{itemize}
\item[$(a)$]
Any definable action of $G$ on a finite set is trivial. 
\item[$(b)$]
If $G'$ is finite, then $G$ is abelian. 
\item[$(c)$] 
If $Z(G)$ is finite, then $G/Z(G)$ is centerless. 
\end{itemize}
\end{fait}
\proof
($a$). The fixator of any element is a definable subgroup of finite index in $G$, and equals $G$ 
by definably connectedness of the latter. 

($b$). Let $g$ be any element in $G$. Notice that the $G$-conjugacy class $g^G$ lies in the coset $gG'$, 
since $G/G'$ is abelian. Since $G'$ is finite, the coset $gG'$ is finite and it follows that the $G$-conjugacy 
class $g^G$ is finite as well. As $G$ acts definably on $g^G$ by conjugation, claim ($a$) shows that $g^G$ is central in $G$, and in particular $g$ is central in $G$. 

($c$). One may argue exactly as in \cite[Lemma 6.1]{MR1321141}. In fact it suffices to apply 
Fact \ref{FactGenAdMap} to elements $x\in Z_2(G)$, in a fashion very similar to what is done in 
Lemma \ref{FactNilpGpsZGInfStrong} below. 
\qed

\begin{fait}\label{ProductConGps}
Let $A$ and $B$ be two definably connected definable subgroups of a definable group $G$, with $A$ 
normalizing $B$. Then $AB$ is a definably connected definable subgroup of $G$. 
\end{fait}
\proof
Clearly $AB$ is a definable subgroup since $A$ normalizes $B$, and thus we only have to check its 
definable connectedness. 

Let $H$ be a definable subgroup of $AB$ of finite index in $AB$. Since $H\cap A$ is a definable subgroup 
of finite index of $A$, it must be $A$ by definable connectedness of $A$. Hence $A\leq H$. 
Similarly, $B\leq H$. We have shown that any definable subgroup $H$ of finite index in $AB$ satisfies 
$AB\leq H$. Hence $AB$ has no proper definable subgroup of finite index, and thus it is definably connected.  
\qed

\bigskip
From now on we consider only groups with the $dcc$. 
If $H$ is a definable subgroup of a group with the $dcc$ on definable subgroups, 
then $H$ has a smallest definable 
subgroup of finite index, the intersection of all of them. It is denoted by $H\o$ and 
called the {\em definably connected component} of $H$. 
Clearly $H\o$ is definably characteristic in $H$, and in particular normal in $H$. 

Any subset $X$ of a group $G$ with the $dcc$ on definable subgroups 
is contained in a smallest definable subgroup, the intersection 
of all of them. It is called the {\em definable hull} of $X$ and denoted by $H(X)$. 

\begin{rem}\label{RemXKinv}
Let $G$ be a definable group with the $dcc$ on definable subgroups, $X$ and $K$ two 
arbitrary subsets of $G$ with $X$ $K$-invariant. Then $H(X)$ is $K$-invariant. 
\end{rem}
\proof
We have $X^k=X$ for any element $k$ of the arbitrary subset $K$. It follows that 
$X\subseteq {H(X)\cap H(X)^k}$ for every such element $k$. Since $H(X)\cap H(X)^k$ is a definable subgroup 
containing $X$, the definition of $H(X)$ implies that $H(X)=H(X)\cap H(X)^k$, and thus $H(X)\leq H(X)^k$. 
If $H(X)<H(X)^k$, then conjugating by $k^{-1}$ we find 
$$X\subseteq H(X)^{k^{-1}}<H(X)$$ 
(since $X=X^{k^{-1}}$), a contradiction to the definition of $H(X)$. Hence 
$H(X)=H(X)^k$. 
\qed

\bigskip
We now restate a typical fact about groups of finite Morley rank which indeed uses only the 
$dcc$ on definable subgroups. 

\begin{fait}\label{NilpResOfDefHulls}
Let $G$ a group with the $dcc$ on definable subgroups. If $X$ and $Y$ are two arbitrary 
subgroups of $G$, then $[H(X),H(Y)]\leq H([X,Y])$. 
In particular if $X$ is an $n$-nilpotent (resp. $n$-solvable) subgroup, then $H(X)$ is $n$-nilpotent 
(resp. $n$-solvable) as well. 
\end{fait}
\proof
Since $X$ and $Y$ normalize $[X,Y]$, they also normalize $H([X,Y])$ by Remark \ref{RemXKinv}. 
As $N(H([X,Y]))$ is definable, we get that $H(X)$ and $H(Y)$ are both in $N(H([X,Y]))$. Now working 
modulo $H([X,Y])$, the proof works formally exactly as in \cite[Lemma 5.37]{MR1321141} 
(using the notation ``$H$" instead of the notation ``$d$" used there for the definable hull). 
\qed

\bigskip
We mention that Fact \ref{NilpResOfDefHulls} was also observed in 
\cite[Lemma 6.8]{MR2006422} when the ambient group $G$ is definable in an o-minimal structure, 
referring to the same proof. 

If $G$ is any group, then $F(G)$ denotes the subgroup generated by all normal nilpotent subgroups 
of $G$ and is called the {\em Fitting} subgroup of $G$. It is clearly a normal subgroup of $G$. 
For the final fact of this section we use a theorem valid for all {\em ${\cal M}_c$-groups}, that is groups 
with the mere $dcc$ on centralizers, as opposed to the $dcc$ on all definable subgroups. 

\begin{fait}\label{faitFGNilpDef}
Let $G$ be a group with the $dcc$ on centralizers. Then $F(G)$ is nilpotent. 
In particular if $G$ is a group with the $dcc$ on definable subgroups, then $F(G)$ is nilpotent and definable. 
\end{fait}
\proof
The nilpotency of $F(G)$ for all ${\cal M}_c$-groups is shown in \cite[Theorem 1.2.11]{MR1473226}. 
When $G$ has the $dcc$ on all definable subgroups, the definable hull $H(F(G))$ exists. Since  
$F(G)\trianglelefteq G$, Remark \ref{RemXKinv} implies that $H(F(G))$ is normal in $G$ as well. 
But by Fact \ref{NilpResOfDefHulls} $H(F(G))$ is also nilpotent. Hence in this case 
$H(F(G))\leq F(G)$, and both groups are equal and definable. 
\qed

\section{Nilpotent groups with $dcc$}\label{NilpGpsdcc}

We continue in the same spirit as in the previous section, analyzing the structure of definable 
groups with the $dcc$, now specializing to the case of nilpotent and abelian groups. As in the previous section, 
definability may refer here to definability from ${\cal M}^{\eq}$ for some fixed structure $\cal M$. 
We start with a lemma 
inspired by \cite[Ex. 5 p. 98]{MR1321141} in the finite Morley rank case. 

\begin{lem}\label{FactNilpGpsZGInfStrong}
Let $G$ be a nilpotent group with the $dcc$ on definable subgroups, and 
$H$ an infinite normal subgroup of $G$. Then $H\cap Z(G)$ is infinite. 
\end{lem}
\proof
We proceed by induction on the nilpotency class of $G$. Of course it starts with $G$ abelian, 
and thus we may assume that our counterexample of minimal nilpotency class is nonabelian. 

We suppose $H\cap Z(G)$ finite. Applying the induction hypothesis in $G/Z(G)$, 
which has smaller nilpotency class by Fact \ref{FactZnGnNilpGps} and clearly also 
satisfies the $dcc$, 
we see that the image of $H$ in $G/Z(G)$ has an infinite intersection with $Z(G/Z(G))$. 
In other words, 
$$HZ(G)/Z(G)\cap Z(G/Z(G))$$ 
is infinite. 
Replacing $H$ by $H\cap Z_2(G)$, we may thus assume $H\leq Z_2(G)$ without loss of generality. 

Let $x\in G$. Then the application $\ad_x : h  \longmapsto [h,x]$, defined for $h\in H$, takes its values in 
$H \cap Z(G)$ since $H\trianglelefteq G$ and $H \leq Z_2(G)$. By Fact \ref{FactGenAdMap} this is 
a group homomorphism, with $C_H(x)$ as kernel. As $H\cap Z(G)$ is finite, its image is finite, and 
thus the kernel $C_H(x)$ has finite index in $H$. 

Besides, $Z(G)=C_G(x_1, \cdots ,x_n)$ for finitely many elements $x_1$,..., $x_n \in G$ by $dcc$ on centralizers. 
In particular 
$${H \cap Z(G)}={C_H(x_1) \cap \cdots \cap C_H(x_n)}$$ 
has finite index in $H$. But since $H$ is infinite and $H\cap Z(G)$ is finite this is a contradiction. 
\qed

\bigskip
It is a question whether Lemma \ref{FactNilpGpsZGInfStrong} can be strengthened to groups with 
the mere $dcc$ on centralizers. We note here that there exist groups $G$ with the $dcc$ on centralizers 
but such that $G/Z(G)$ does not have this property \cite{MR549936}. 

We now prove an ``infinite" version of the classical normalizer condition for finite nilpotent groups. 
We will not use it in the present paper but we expect it to be as useful as in  \cite{MR2225896} for a 
potential theory of Carter subgroups in groups definable in o-minimal structures. 

\begin{lem}
Let $G$ be a nilpotent group with the $dcc$ on definable subgroups, and 
$H$ a definable subgroup of infinite index in $G$. Then $N_G(H)/H$ is infinite. 
\end{lem}
\proof
Our proof can be compared to that of \cite[Lemma 6.3]{MR1321141} in the finite Morley rank case, but here 
we argue by induction on the nilpotency class rather than on the dimension, and we deal with 
``non-connected finite bits" throughout. We note also that we do not need the full $dcc$ on definable 
subgroups, but merely the existence of definably connected components of definable subgroups. 

Our counterexample $G$ of minimal nilpotency class is of course not abelian. 
Take such a $G$ and  the corresponding $H\leq G$.  Let $Z=Z(G)$. 
Since $G$ is a counterexample to our statement and since  $Z\o\leq N(H)$, we have 
$Z\o\leq H\o \leq H$. If $HZ$ had finite index in $G$, then we would get $H\o=(HZ)\o=G\o$, and 
then $H$ could not be of infinite index in $G$.  This shows that $HZ$ has infinite index in $G$, and it 
follows that $HZ/Z$ has infinite index in $G/Z$ as well. 
By Fact \ref{FactZnGnNilpGps}, $G/Z$ has smaller nilpotency class, and the induction hypothesis applied 
in this quotient now implies that $HZ$ has infinite index in $N(HZ)$ (since $N(HZ)$ is exactly the preimage 
in $G$ of $N_{G/Z}(HZ/Z)$). It follows that $(HZ)\o$ has infinite index in $N\o(HZ)$ and since 
$H\o=(HZ)\o$ we get 
$$[N\o(HZ):H\o]=\infty.$$
We are done if we show that $N\o(HZ)$ normalizes $H$. Of course, $N\o(HZ)$ normalizes $(HZ)\o=H\o$. 
Take now any element $h$ in $H$. Then the coset $hH\o$ is also a coset of the form $h(HZ)\o$, with $h$ in 
$HZ$. By Fact \ref{ActionConOnFinite} $(a)$, the definably connected group 
$N\o(HZ)$ induces by conjugation a trivial action on the finite quotient $(HZ)/(HZ)\o$. In other words, 
$N\o(HZ)$ setwise stabilizes by conjugation the coset $hH\o$. We have shown that $N\o(HZ)$ normalizes 
$H$, and this finishes our proof. 
\qed

\bigskip
We now pass to a finer analysis of abelian groups with the $dcc$, essentially as in 
\cite{MR0289280} in the finite Morley rank case; see also \cite[Theorem 6.7]{MR1321141}. 

\begin{fait}\label{MacintyresThmAbGps}
Let $A$ be an abelian group with the $dcc$ on definable subgroups. Then 
\begin{itemize}
\item[$(a)$]
$A=D\oplus B$ for some definable divisible subgroup $D$ and some (interpretable but 
not necessarily definable) complement $B$ of bounded exponent. 
\item[$(b)$]
$A=DC$ for some definable and characteristic subgroups $D$ and $C$, with $D$ divisible as 
above and $C$ of bounded exponent. Furthermore $D\cap C$ is finite whenever $A$ contains 
no infinite elementary abelian $p$-subgroup. 
\end{itemize}
\end{fait}
\proof
($a$). Let $D=\bigcap_{n\in \N}A(n!)$, where $A(k)=\{a^k~|~a\in A\}$. By $dcc$, $D=A(n!)$ for some 
$n$ and $D$ is definable. 
It is clearly divisible. Now, since $A$ is abelian, a well known theorem of 
Baer implies that any subgroup disjoint from the divisible subgroup $D$ can be extended 
to a complement. In particular 
$A=D\oplus B$ for some complement $B$. Now, clearly, $B$ must be of bounded exponent. 

($b$). Starting from the decomposition $A=D\oplus B$ as in claim ($a$), we can now take 
$C=\{a\in A~|~a^{\exp(B)}=1\}$. Then all our statements are clear, knowing the decomposition 
of any abelian divisible group $D$ as a direct product of quasicyclic Pr\"{u}fer $p$-groups 
and of copies of $\Q$. 
\qed

\begin{cor}\label{ADefCondccNoInfpElemAdiv}
Let $A$ be a definably connected abelian group with the $dcc$ on definable subgroups, and 
with no infinite elementary abelian $p$-subgroup. Then $A$ is divisible. 
\end{cor}
\proof
Let $A=DC$ be a decomposition of $A$ as in Fact \ref{MacintyresThmAbGps} ($b$), 
with $D$ divisible and $C$ of bounded exponent. For $p$ any prime and any $n\geq 0$, let 
$C_{p^n}={\{c\in C~|~g^{p^n}=1\}}$. Since $C$ has bounded exponent, the $p$-primary component 
of $C$ has the form $C_{p^N}$ for some $N$, and it is in particular definable. 
By assumption, $C_p$ is finite. 
Considering the group homomorphism $x \longmapsto x^p$, we see by
induction on $n$ that each subgroup $C_{p^n}$ is finite. 
In particular, the $p$-primary component $C_{p^N}$ of $C$ is finite. 

Suppose now towards a contradiction $D<A$. Since $A=DC$ and $A$ is definably connected, we then 
get that $C$ is infinite. Hence $C$ is the direct product, for primes $p$ varying in an infinite set $I$, 
of nontrivial finite $p$-groups, each of exponent $p^N$ for some $N$ depending on $p$. Now taking 
successive $p^N$th powers of $C$, as $p$ varies in $I$, we build an infinite descending chain of definable 
subgroups of $C$, as contradiction. Hence $A=D$ is divisible. 
\qed

\bigskip
As a consequence of Fact \ref{MacintyresThmAbGps} we also get a result of lifting of torsion. It can 
be compared to \cite[Ex. 11 p. 93]{MR1321141}, which have had endless applications 
in the finite Morley rank case. 

\begin{fait}\label{FactLiftTorsion}
Let $G$ be a group with the $dcc$ on definable subgroups, 
$N$ a normal definable subgroup of $G$, and $x$ an element of $G$ such that $x$ 
has finite order $n$ modulo $N$. Then the coset $xN$ 
contains an element of finite order, involving the same prime divisors as in $n$.  
\end{fait}
\proof
The definable hull $H(x)$ is abelian by Fact \ref{NilpResOfDefHulls}. 
Replacing $G$ by $H(x)$ and $N$ by $N\cap H(x)$, we may assume $G$ abelian. 
Now from Fact \ref{MacintyresThmAbGps} we get a decomposition 
$N=D\oplus B$ where $D$ is $n$-divisible and $B$ is a direct sum of $p$-groups, for $p$ prime 
dividing $n$ (we may ``transfer" from the original $B$ to the original $D$ given by Fact \ref{MacintyresThmAbGps} ($a$) the Sylow $q$-subgroups for $q$ not dividing $n$). 
Since $x^n$ is in $N$, $x^n=d^nb$ for some $d$ in $D$ and some $b$ in $B$. Now 
$xd^{-1}$ is in the same $N$-coset as $x$, and since $(xd^{-1})^n=b$ the element $xd^{-1}$ satisfies 
our requirement. 
\qed

\bigskip
By \cite[Theorem 5.1]{MR2436142} (see also \cite[Proposition 6.1]{MR1303545}), 
any abelian group of bounded exponent {\em definable} in an o-minimal structure must be finite. 
Hence any definably connected abelian group definable in an o-minimal structure is divisible, 
by Fact \ref{MacintyresThmAbGps} ($b$) or Corollary \ref{ADefCondccNoInfpElemAdiv} 
(contrarily to the finite Morley rank case where infinite bounded exponent subgroups are a real 
possibility). 

In the finite Morley rank case, Fact \ref{MacintyresThmAbGps} ($b$) generalizes identically 
to nilpotent groups by a theorem of Nesin. For nilpotent groups $G$ definable in o-minimal structures, 
we always have that $G\o$ is divisible by \cite[Theorem 6.10]{MR2006422}. In any case, all 
assumptions made in the next lemma, and in the reminder of this section, are met for groups 
definable in o-minimal structures.

\begin{lem}\label{StructGenGpsNilp}
Let $G$ be a nilpotent group with the $dcc$ on definable subgroups and with $G\o$ divisible. Let 
$\Tor(G)$ denote the set of torsion elements of $G$. Then
\begin{itemize}
\item[$(a)$] 
$\Tor(G)$ is a subgroup of $G$, and it is the direct product of its Sylow $p$-subgroups. 
Furthermore it commutes with $G\o$. 
\item[$(a')$] 
$\Tor(G\o)$ is in $Z\o(G\o)$, and in particular it is a divisible abelian subgroup of $G\o$.  
Besides, any definable section of $G\o/Z\o(G\o)$ is torsion-free. 
\item[$(b)$] 
$G=G\o \ast B$ (central product) for some finite subgroup $B$ of $\Tor(G)$. 
\item[$(c)$] 
If $G\o$ has no infinite elementary abelian $p$-subgroup, then the subgroup $B$ as in claim 
(b) may be chosen to be a finite characteristic subgroup. 
\end{itemize}
\end{lem}
\proof
($a$) The fact that $\Tor(G)$ is a subgroup and the direct product of its Sylow $p$-subgroups is 
true in any nilpotent group (\cite[5.2.7]{MR1357169}). 

Since $G\o$ is $p$-divisible for any prime $p$, 
it commutes with any Sylow $p$-subgroup by \cite[4.2; 4.7; 6.13]{MR0409661}. 
In particular $\Tor(G)$ commutes with $G\o$. 

($a'$) In particular, the set of torsion elements of $G\o$ forms a central subset of $G\o$, 
and hence it is a central subgroup of $G\o$. By divisibility of $G\o$, 
$\Tor(G\o)$ is also a divisible (abelian) subgroup of $G\o$. 
Hence any finite quotient of $\Tor(G\o)$ must be trivial, 
and thus $\Tor(G\o)Z\o(G\o)/Z\o(G\o)$ is trivial. Hence $\Tor(G\o)\leq Z\o(G\o)$. 
Our last claim follows then from the lifting of torsion given by Fact \ref{FactLiftTorsion}. 

($b$) Any extension of a locally finite group by a locally finite group 
is locally finite \cite[Lemma 1.A.2 p. 2]{MR0470081}, 
and in particular any torsion solvable group is locally finite. 
It follows that $\Tor(G)$ is locally finite. But by lifting of torsion, Fact \ref{FactLiftTorsion}, we get that 
$G=G\o\ast \Tor(G)$, since $G/G\o$ is finite. Now we can also pick up finitely many 
representatives of finite orders of all cosets of $G\o$ in $G$; they generate a finite subgroup 
$B$ of $\Tor(G)$ by local finiteness of the latter, and we also have 
that $G=G\o \ast B$. 

($c$) 
Let $n$ be the order of $B$, or the minimum common multiple of orders of elements of $B$. 
Consider the set $\tilde B$ of elements $g$ in $\Tor(G)$ such that $g^n=1$. Clearly $B\subseteq \tilde B$. 
Any element $g$ in $\tilde B$ has the form $bh$ for some $b\in B$ and some element $h$ in $G\o$. 
Since $[b,h]=1$, we see that $1=g^n=(bh)^n=b^nh^n=h^n$. 

Now the torsion subgroup of $G\o$ consists of direct products of quasicyclic Pr\"{u}fer $p$-groups. 
Our extra assumption that $G\o$ contains no infinite elementary abelian $p$-subgroup for any prime $p$ 
now implies that each $p$-primary component of $G\o$ is the direct product of at most finitely many 
copies of the quasicyclic Pr\"{u}fer $p$-group. 
In particular, for $n$ as above, there are only finitely many elements $h$ in $G\o$ satisfying $h^n=1$. 
Since $B$ is finite, this shows that $\tilde B$ is also finite. 

Since $\tilde B$ is a finite subset of $\Tor(G)$, we conclude by local finiteness of the latter that 
$\<\tilde B\>$ is finite. Clearly, by definition, $\<\tilde B\>$ is characteristic in $G$. So we may replace 
$B$ by $\<\tilde B \>$. 
\qed

\begin{cor}\label{CorHNormalNilpHoZH}
Let $G$ be a definably connected group with the $dcc$ on definable subgroups. Let 
$H$ be a definable normal nilpotent subgroup of $G$, with $H\o$ divisible and without 
infinite elementary abelian $p$-subgroups. Then $H=(Z(G)\cap H)H\o$, and $\Tor(H) \leq Z(G)$. 
\end{cor}
\proof
The first point follows from Lemma \ref{StructGenGpsNilp} ($c$) and Fact \ref{ActionConOnFinite} ($a$). 

For the second point, 
we have seen in the proof of Lemma \ref{StructGenGpsNilp} ($c$) that the set of elements of $H\o$ 
of order dividing $n$ is finite, for every $n$, and thus Fact \ref{ActionConOnFinite} ($a$) also gives that the 
torsion of $H\o$ is central in $G$. Now the full torsion subgroup $\Tor(H)$ of $H$ is in $Z(G)$ 
by Lemma \ref{StructGenGpsNilp} ($c$) and Fact \ref{ActionConOnFinite} ($a$) as well. 
\qed

\bigskip
Of course, Corollary \ref{CorHNormalNilpHoZH} may apply with $H$ the Fitting subgroup of $G$. 

\begin{cor}
Let $G$ be a definably connected group with the $dcc$ on definable subgroups. 
Suppose that $F\o(G)$, which is definable and nilpotent by Fact \ref{faitFGNilpDef}, 
is divisible and without infinite elementary abelian $p$-subgroups. 
Then $F(G)=Z(G)F\o(G)$, and $\Tor(F(G)) \leq Z(G)$. 
\end{cor}
\proof
Since $Z(G)\leq F(G)$, Corollary \ref{CorHNormalNilpHoZH} applies directly. 
\qed

We conclude the present section with a flexible characterization of central extensions 
in the general context under consideration. Again, for the reasons mentioned before 
Lemma \ref{StructGenGpsNilp}, everything here applies to any group definable in an o-minimal structure. 

\begin{lem}\label{FactCentralExtNonCon}
Let $G$ be a definably connected group with the $dcc$ on definable subgroups, 
and let $H$ be a definable normal subgroup such that $H\o$ contains no infinite elementary 
abelian $p$-subgroups. Then $H\leq Z(G)$ if and only if $H\o\leq Z(G)$. 
\end{lem}
\proof
Suppose $H\o\leq Z(G)$. 
Since $H$ is finite modulo $H\o$, Fact \ref{ActionConOnFinite} ($a$) 
implies that $H/H\o \leq Z(G/H\o)$. Hence $[G,H]\leq H\o$, and in particular 
$H'\leq H\o$. 
But by assumption $H\o \leq Z(G)$, so in particular we get $H'\leq Z(H)$. So the definable 
subgroup $H$ is $2$-nilpotent. Besides, the definably connected abelian subgroup $H\o$ is 
divisible by Corollary \ref{ADefCondccNoInfpElemAdiv}. 
Now by Corollary \ref{CorHNormalNilpHoZH} we get $H\leq Z(G)$. 
\qed

\section{Groups with a definable and additive dimension}\label{SectionDefAddDimension}

We continue as in the two previous sections with rather general considerations. We will see in the next 
section that this applies to groups definable in o-minimal structures. Throughout the present section, 
$G$ is a group interpretable in a structure $\cal M$, and definability refers to ${\cal M}^{\eq}$. 
We suppose that to each definable set in Cartesian products of $G$ is attached a dimension in $\N$ 
and denoted by $\dim$. We essentially require 
the dimension to be definable and additive, and also require a natural additional 
axiom on the dimension of finite sets. Axioms are the following. 
\begin{itemize}
\item[(A1)]{\bf (Definability)} 
If $f$ is a definable function between two definable sets $A$ and $B$, then for every $m$ in $\N$ 
the set $\{ b\in B~|~\dim(f^{-1}(b))=m \}$ is a definable subet of $B$. 
\item[(A2)]{\bf (Additivity)} 
If $f$ is a definable function between two definable sets $A$ and $B$, whose fibers have constant 
dimension $m$ in $\N$, then $\dim(A)=\dim(B)+m$. 
\item[(A3)]{\bf (Finite sets)}
A definable set $A$ is finite if and only if $\dim(A)=0$.
\end{itemize}

Axioms A2 and A3 guarantee that if $f$ is a definable bijection between two definable sets $A$ and $B$, 
then $\dim(A)=\dim(B)$. In other words, they guarantee that the dimension is preserved under 
definable bijections. 

Using axiom A2, one sees then that if $K\leq H\leq G$ are definable subgroups of $G$, then 
$$\dim(H)=\dim(K)+\dim(H/K),$$ 
a special case of the additivity we will use freely throughout. 
 
In general axioms A1-2 suffice for computing the dimensions of unions of uniformly definable 
families of sets (see \cite[Fact 4 and Corollary 5]{MR2537672} in the finite Morley context). 
Here we will merely use a computation of the dimension for products of definable subgroups. 

\begin{lem}\label{LemDimProduct}
Let $G$ be a group equipped with a dimension satisfying axioms A1-3. 
If $A$ and $B$ are two definable subgroups of $G$, then 
$\dim(AB)=\dim(A)+\dim(B)-\dim(A\cap B)$. 
\end{lem}
\proof
It suffices to consider the definable map $(a,b)\longmapsto ab$ from $A\times B$ onto $AB$. 
Its fibers are in definable bijection with $A\cap B$. 
\qed

Clearly any group equipped with a dimension satisfying axioms A1-3, and which satisfies the $dcc$ on definable subgroups of the same dimension, satisfies the $dcc$ on all definable subgroups. 
Indeed, if $H_1>H_2>\cdots$ 
is a strictly descending chain of definable subgroups, then subchains of subgroups with the same dimension 
must be finite by assumption, and chains of subgroups with strictly decreasing dimensions must be finite as well. 

In the general context of groups with a dimension satisfying axioms A1-3 
we do not have a good version of Zilber's stabilizer argument, and consequently Zilber's generation lemma 
for indecomposable sets as in the finite Morley rank context (see \cite[Theorem 5.26]{MR1321141}). 
Our work here will consist mainly in reductions to the following lemma, which can be seen as 
a rudimentary version of Zilber's generation lemma. 

\begin{lem}\label{WeakZilberIndec}
Let $G$ be a group equipped with a dimension satisfying axioms A1-3, and let 
$A_i$, $i\in I$, be an arbitrary family of definably connected definable subgroups of $G$ which 
pairwise normalize each others. 
Then $\<A_i~|~i\in I\>=A_{i_1}\cdots A_{i_k}$ for finitely many $i_1$, ...$i_k \in I$, and in particular it is definable. Furthermore we can always take $k$ to be smaller than or equal to the dimension of that definable subgroup. 
\end{lem}
\proof
First notice that the definable subgroups $A_i$ of $G$ can be definably connected even if $G$ does not 
satisfy the $dcc$, as in the context of Fact \ref{ProductConGps}. 

Any finite product $A_{i_1}\cdots A_{i_k}$ is a definably connected and definable subgroup by 
Fact \ref{ProductConGps}. Now one sees with Lemma \ref{LemDimProduct} 
that any such finite product $A_{i_1}\cdots A_{i_k}$ 
of maximal dimension is equal to $\<A_i~|~i\in I\>$. 

Once we know that $\<A_i~|~i\in I\>=A_{i_1}\cdots A_{i_k}$ is definable (and definably connected), we see 
again with Lemma \ref{LemDimProduct} that we may choose $k\leq \dim(\<A_i~|~i\in I\>)$. 
\qed

\begin{question}
Can one relax the normalization hypothesis in Lemma \ref{WeakZilberIndec}, assuming the ambient 
group, say, nilpotent? 
\end{question}

\bigskip
The following corollary of Lemma \ref{WeakZilberIndec} is going to be the basic tool for the proof of 
our main theorem. 

\begin{cor}\label{Lem[X,H]definable}
Let $G$ be a group equipped with a dimension satisfying axioms A1-3, 
$H$ a definably connected definable subgroup, and $X$ an 
arbitrary subset of $G$. Suppose that $[H,X]$ is abelian and centralized by $H$. 
Then $[H,X]$ is a definably connected definable (abelian) subgroup of $C(H)$, and actually 
we have $[H,X]=[H,X]_{\dim([H,X])}$
\end{cor}
\proof
Fix any element $x$ in $X$. Then the map $\ad_x~:~H \longrightarrow G$ is a group homomorphism 
by Fact \ref{FactGenAdMap}, and clearly it is definable. 
Hence its image $[H,x]$ is a definably connected definable subgroup of $G$. 

Now all such subgroups $[H,x]$ pairwise normalize each other  since they all live in 
$[H,X]$ which is abelian. But $[H,X]$ is generated by these subgroups $[H,x]$, and hence it suffices 
now to apply Lemma \ref{WeakZilberIndec}. 
It implies that $[H,X]$ is a definably connected definable (abelian) subgroup of $G$. 

Our last remark about the commutator width of $[H,X]$ corresponds to the last claim in 
Lemma \ref{WeakZilberIndec}. 
\qed

If $G$ is any group, $R(G)$ denotes the subgroup generated by all normal solvable subgroups of $G$ and 
is called the {\em solvable radical} of $G$. We start with a general remark about solvable radicals 
in arbitrary definably connected groups. 

\begin{rem}\label{RemRGFinite}
Let $G$ be any definably connected group. 
If $R(G)$ is finite, then $R(G)=Z(G)$ and $G/R(G)$ has a trivial solvable radical. 
\end{rem}
\proof
The first point follows from Fact \ref{ActionConOnFinite} ($a$) and the second is obvious once we know 
that $R(G)$ is solvable. 
\qed

We now prove the definability and the solvability of the solvable radical in a rather abstract context, which 
incorporates both groups of finite Morley rank and groups definable in o-minimal structures. 

\begin{lem}\label{DefSolvofRG}
Let $G$ be a group equipped with a dimension satisfying axioms A1-3 and with the $dcc$ on definable 
subgroups of same dimension. Then $R(G)$ is definable and solvable. 
\end{lem}
\proof
By assumption $G$ satisfies the $dcc$ on all definable subgroups. 

By Remark \ref{RemXKinv} and Fact \ref{NilpResOfDefHulls}, 
$R(G)$ is generated by the set of {\em definable} normal solvable 
subgroups of $G$, so it suffices to show that the group generated by such definable normal solvable 
subgroups is solvable. By Lemma \ref{WeakZilberIndec}, the 
product of all {\em definably connected} definable normal solvable subgroups 
of $G$ is the product of finitely many 
of them, and hence it is definable and solvable. Let $K$ denote this maximal definably connected definable 
normal solvable subgroup of $G$ (which, of course, is going to be $R\o(G)$). 

Any definable normal solvable subgroup $N$ of $G$ has a finite image in $G/K$. 
So replacing $G$ by $G/K$ we may assume that all normal solvable subgroups are finite. 
If $N$ is such a normal solvable subgroup, then as $N$ is 
normalized by $G\o$ we get that $N\leq C_G(G\o)$ by Fact \ref{ActionConOnFinite} ($a$). 
But $C_G(G\o)$ is a finite subgroup of $G$, as otherwise 
$C\o_G(G\o)$ is a nontrivial definably connected definable subgroup of $G$, 
hence in $G\o$, and then in the center of $G\o$ which is finite, a contradiction. 
Hence any normal solvable subgroup of 
$G$ is in the finite (normal) subgroup $C_G(G\o)$ of $G$, and in particular in $R(C_G(G\o))$, which is solvable. 
This completes the proof. 
\qed

\begin{rem}
Lemma \ref{DefSolvofRG} generalizes the definability and the solvability of $R(G)$, 
noticed by Belegradek in the finite Morley rank context. As we will see in the next section 
that groups definable in o-minimal structures satisfy all our current assumptions, 
it also gives the definability and the solvability of the 
full solvable radical $R(G)$ in this context. We note that in the literature on groups in o-minimal structures it 
is usually only $R\o(G)$ which is defined, essentially 
by the same argument as in the second paragraph of 
the proof of Lemma \ref{DefSolvofRG}. See for example \cite[Fact 2.3]{ConversannoArticle1}. 
\end{rem}

\bigskip
We finish the present section with a mere definition. If $G$ is a group as in 
Lemma \ref{DefSolvofRG}, then $R(G)$ is a solvable definable subgroup, and $F(G)$ is a 
nilpotent definable subgroup by Fact \ref{faitFGNilpDef}. 
Of course, we have 
$$F(G)\trianglelefteq R(G)\trianglelefteq G$$
and the same inclusions hold with definably connected components. If $G$ is definably connected 
and $R(G)$ is finite (equivalently $R\o(G)=1$), then Remark \ref{RemRGFinite} shows that 
$Z(G)=F(G)=R(G)$ and $G/R(G)$ has no nontrivial normal abelian subgroup. In this 
case $G$ can have only finite normal abelian subgroups (in $Z(G)$). 

The definition of semisimplicity for groups may vary in the literature, depending on whether authors 
admit a finite center or not. Here we will admit such finite centers in our definition of semisimplicity. 

\begin{defi}
A definably connected group $G$ as in Lemma \ref{DefSolvofRG} is said to be 
{\em semisimple} whenever $R\o(G)=1$. By Remark \ref{RemRGFinite}, this is equivalent to requiring 
that $Z(G)=R(G)$ is finite and $G/R(G)$ has trivial abelian normal subgroups. 
\end{defi}

\section{Groups definable in o-minimal structures}\label{Section6}

From now on we specialize our analysis to groups {\em definable} in o-minimal structures. Here we insist 
on definability, as opposed to interpretability, because since \cite{MR961362} 
the whole theory of groups in o-minimal structures 
has been developed in this restricted context. Hence the framework from now on 
is the following: $\cal M$ is an o-minimal structure and $G$ is a group definable in $\cal M$. In general 
the o-minimal structure $\cal M$ itself may not eliminate imaginaries, but in some sense $G$ does. In fact, 
$G$ has {\em strong definable choice} in the following sense. 

\begin{fait}\label{FactDefiChoice}
{\bf \cite[Theorem 7.2]{MR2006422}} 
Let $G$ be a group definable in an o-minimal structure $\cal M$, and let $\{T(x)~|~x\in X\}$ be a 
definable family of non-empty definable subsets of $G$. Then there is a definable function 
$t~:~X\longrightarrow G$ such that for all $x$, $y\in X$ we have $t(x)\in T(x)$ and if 
$T(x)=T(y)$ then $t(x)=t(y)$. 
\end{fait}

\bigskip
In particular, a definable section $H/K$ of a group $G$ definable in an o-minimal structure $\cal M$ 
can also be regarded as definable in $\cal M$, even though it has a priori to be considered as definable 
in ${\cal M}^{\eq}$. Hence when making proofs by induction below we will freely pass to definable sections 
of groups definable in o-minimal structures, and still consider them as definable. 

As seen already in Sections \ref{Gpsdcc} and \ref{NilpGpsdcc}, groups definable in o-minimal structures 
satisfy the $dcc$ on definable subgroups, and all the results of these two sections apply to them. 
Similarly, all results of Section \ref{SectionDefAddDimension} apply to groups definable 
in o-minimal structures: for the definition of the dimension of definable subsets, which uses the cell decomposition, we refer to \cite{MR1633348}. 
The definability and the additivity (our axioms A1 and A2 respectively) 
can be found in \cite[Chapter 4 (1.5)]{MR1633348}, and the fact that finite sets are 
exactly the 0-dimensional sets (our axiom A3) is noticed in \cite[Chapter 4 (1.1)]{MR1633348}. 

A group is {\em definably simple} if it is nonabelian and has no proper nontrivial normal definable 
subgroup (in the specific universe considered, typically our ground o-minimal structure here). 
As seen at the end of Section \ref{SectionDefAddDimension}, if $G$ is a group definable in an o-minimal 
structure then $R(G)$ is definable and solvable. When $G$ is definably connected, then 
$G/R\o(G)$ is semisimple and $G/R(G)$ has a nice description, following a solution 
of a version of the Cherlin-Zilber conjecture for groups in o-minimal structures. 

\begin{fait}\label{FactStrucSemisimpleGps}
{\bf \cite[Theorem 4.1]{MR1707202}}
Let $G$ be an infinite definably connected group definable in an o-minimal structure, with $R(G)=1$. 
Then $G$ is the direct product of    definable  definably simple subgroups $G_1$, ..., $G_k$. 
\end{fait}

For Fact \ref{FactStrucSemisimpleGps}, see also \cite[Remark 4.2]{MR1707202}. 
There is in fact much more information known about the (finitely many) definably simple 
factors $G_i$ appearing in Fact \ref{FactStrucSemisimpleGps}, and we can see them in 
at least two ways. On the one hand, for any definably simple group $G$ definable in an o-minimal 
structure $\cal M$, there is a definable real closed field $\cal R$ and a real algebraic group $H$ 
defined over the subfield of real algebraic numbers $R_{alg}\subseteq R$, such that $G$ 
is definably isomorphic in $\cal M$ to $H(R)\o$, the definably connected component of $H(R)$ 
(See \cite[4.1]{MR1707202} for the existence of $H$ and the proof of 
\cite[5.1]{MR1873380} for the fact that $H$ can be defined over $R_{alg}$; 
see also \cite[Fact 1.1]{HPPCentExt}). On the other hand, we can see $G$ as a group elementarily 
equivalent to a simple (centerless) Lie group by \cite[Theorem 5.1]{MR1873380}.  

We now prove that such definably simple groups have finite commutator width. Recall that a 
group $G$ is {\em perfect} if $G'=G$. 

\begin{fait}\label{FactGDefSimpleGSimple}
Let $G$ be a definably simple group definable in an o-minimal 
structure. Then  $G=[G,G]_k$ for some $k$, and in particular $G=G'$ is perfect. 
\end{fait}
\proof
Clearly, it suffices to prove the first statement. 

If $G$ is finite, then of course the definable simplicity of $G$ implies its abstract group theoretic 
simplicity. Since $G'$ is a nontrivial normal subgroup, we get $G=G'=[G,G]_k$ for some $k$. 

Suppose now $G$ infinite. Our assumptions pass to elementary extensions, so taking a sufficiently 
saturated elementary extension of the ground o-minimal structure, we can work in a 
sufficiently saturated elementary extension $G^*$ of $G$. If we show there that ${G^*}'=[G^*,G^*]_k$ 
for some $k$, then ${G^*}'$ is definable, and $G^*={G^*}'=[G^*,G^*]_k$ by definable simplicity of $G^*$, 
showing then by elementary equivalence that $G=[G,G]_k$ as well. 
So we may assume that $G$ is already sufficiently saturated 
(actually $\omega$-saturated suffices in what follows). 

Consider first the case $G$ not definably compact. 
Then by \cite[Corollary 6.3]{MR1873380} $G$ is abstractly simple as a group. 
Since $G'$ is a nontrivial normal subgroup of $G$, we then get $G=G'$ by abstract simplicity of $G$. 
Now $G$ is the countable union of the definable sets $[G,G]_k$, so the $\omega$-saturation of $G$ implies 
that $G=[G,G]_k$ for some $k$. Hence in this case we are done. 

Consider now the case $G$ definably compact. As commented after 
Fact \ref{FactStrucSemisimpleGps}, $G$ is  definably isomorphic to a semialgebraic 
group over a real closed field $\cal R$ which is definable in our original o-minimal structure 
(but may be different from it). Hence we can see $G$ as a definably compact group 
definable over an o-minimal expansion of an ordered field, and thus in particular 
over an o-minimal expansion of an ordered group. 
Then we can get $G'=[G,G]_k$ for some $k$ as in the proof of 
\cite[Corollary 6.4 (i)]{HPPCentExt}, and in this case we are done also. 
\qed

\begin{question}
\label{QuestionComWidthDefSimplGp}
Let $G$ be a definably simple group definable in an o-minimal structure. Then is it the case that 
$[G,G]=[G,G]_{\dim(G)}$?
\end{question}

We strongly believe in a positive answer to Question \ref{QuestionComWidthDefSimplGp} but 
unfortunately we haven't found any model-theoretic proof. Actually, it seems even expected that 
the commutator width should be $1$ in Question \ref{QuestionComWidthDefSimplGp}, 
or in other words that every element is a commutator. 
This corresponds to a conjecture of Ore in the case of {\em finite} 
simple groups (proved recently, by cases inspection \cite{OreFiniteSimpleGps}). 
As commented after Fact \ref{FactStrucSemisimpleGps}, 
there are two ways to see a definably simple group as in Question \ref{QuestionComWidthDefSimplGp}, 
either as the semialgebraic connected component of a real algebraic group or as 
(elementarily equivalent to) a simple Lie group. 
When $G$ is definably compact, then it is known that every element is a commutator by an older theorem 
of Got\^o about semisimple Lie groups \cite{MR0033829} 
(see also \cite[Theorem 6.55]{MR2261490}). 
Hence the question reduces to abstractly simple groups, 
as seen in the proof of Fact \ref{FactGDefSimpleGSimple}. 
In general it remains open but it is explicitly conjectured, from the simple Lie groups point of view, 
that $G=[G,G]_1$ in \cite[Conjecture A]{MR836712}. See also \cite[Page 15]{MR1635671} 
for the same question and links with the Cartan decomposition. 
The reader can also consult \cite{MR1422600} for a proof of the conjecture of Ore for groups 
of algebraic or Chevalley type. 

We finish this section with another general remark about the structure of a group definable in 
an o-minimal structure. It is not used in the rest of the paper, but it is certainly 
worth to mention it at this point. 

\begin{rem}
Let $G$ be any group definable in an o-minimal structure. By lifting of torsion, 
Fact \ref{FactLiftTorsion}, any definable subgroup without torsion must be definably connected, and 
the product of two such groups $A$ and $B$ which normalize each other  is also torsion-free: 
if $AB$ contains a torsion element, then the same occurs in $AB/A \simeq{B/(A\cap B)}$. Now 
Lemma \ref{WeakZilberIndec} shows the existence of a unique maximal definable normal torsion 
free subgroup of $G$, the group generated by all of them. 
This gives a somewhat conceptual proof of the existence of such a subgroup, 
noticed recently in \cite{ConversannoArticle1} using more specific machinery 
involving Euler characteristic. 
\end{rem}

\section{Definably connected solvable groups}\label{Section7}

In the present section we consider the structure of {\em solvable} definably connected groups 
definable in o-minimal structures. We first note that there is an o-minimal version of the 
Lie-Kolchin-Mal'cev theorem. 

\begin{fait}\label{GConResG'nilp}
{\bf \cite[Theorem 6.9]{MR2006422}} 
Let $G$ be a definably connected solvable group definable in an o-minimal structure. 
Then $G'$ is nilpotent. 
\end{fait}

As in \cite{MR992314} in the finite Morley rank case, the proof 
of Fact \ref{GConResG'nilp} is based ultimately on linearization. 
It is a question whether a version of Fact \ref{GConResG'nilp} 
can be proved under more general assumptions such as those used in 
Sections \ref{Gpsdcc}, \ref{NilpGpsdcc}, or \ref{SectionDefAddDimension}. 
We also note that our main theorem implies that $G'$ is definable in 
Fact \ref{GConResG'nilp} (Corollary \ref{CorGoSolvGnGndef} below), 
but that Fact \ref{GConResG'nilp} was proved before knowing this. 

Recall that the Fitting subgroup is nilpotent and definable for any group definable in an 
o-minimal structure by Fact \ref{faitFGNilpDef}. 

\begin{prop}\label{GConResGFGDivAb}
Let $G$ be a definably connected solvable group definable in an o-minimal structure. 
Then $G'\leq F\o(G)$. In particular $G/F\o(G)$ and $G/F(G)$ are divisible abelian groups. 
\end{prop}
\proof
Notice that both quotients $G/F\o(G)$ and $G/F(G)$ can then be considered as definable by Fact \ref{FactDefiChoice}, and are both definably connected as $G$ is. 

By Fact \ref{GConResG'nilp}, $G'\leq F(G)$. It follows that 
$G'F\o(G)/F\o(G)$ is finite. 
Now $(G/F\o(G))'=G'F\o(G)/F\o(G)$ is finite, and then it follows 
from Fact \ref{ActionConOnFinite} ($b$) that the definably connected group $G/F\o(G)$ must be abelian. 
Now it follows that $G'\leq F\o(G)$, and $G/F\o(G)$ is abelian, as well as $G/F(G)$. 

The divisibility of $G/F\o(G)$ and $G/F(G)$ is true for any definably connected 
abelian group definable in an o-minimal structure, as mentioned before Lemma \ref{StructGenGpsNilp}. 
\qed

\begin{cor}\label{GConSolvFGInf}
Let $G$ be a nontrivial definably connected solvable group definable in an o-minimal structure. 
Then $F\o(G)$ is nontrivial. In particular $G$ has an infinite abelian characteristic definable subgroup. 
\end{cor}
\proof
Our statement is clear if $G$ is abelian, so suppose $G$ nonabelian. 
Fact \ref{ActionConOnFinite} ($b$) implies that $G'$ is then infinite. 
Now Proposition \ref{GConResGFGDivAb} shows that $F\o(G)$ is infinite as well. For the last 
claim we can use the fact that $F(G)$ has an infinite center, which follows 
from Fact \ref{faitFGNilpDef} and Lemma \ref{FactNilpGpsZGInfStrong}. 
\qed

\bigskip
We are not going to use the following lemma, but it is very similar to an analog fact 
about connected solvable groups of finite Morley rank and which is crucial for a finer 
inductive analysis of such groups. We expect it should have similar consequences for 
definably connected solvable groups definable in o-minimal structures. 

If $H$ and $G$ are two subgroups of a group with $G$ normalizing $H$, then a 
{\em $G$-minimal} subgroup of $H$ is an infinite $G$-normal definable subgroup of $H$, 
which is minimal with respect to these properties (and where definability refers to a fixed 
underlying structure). If $H$ is definable and satisfies the $dcc$ on definable subgroups, then 
$G$-minimal subgroups of $H$ always exist. As the definably 
connected component of a definable subgroup is a definably characteristic subgroup, 
we get also in this case that any $G$-minimal subgroup of $H$ should be definably connected. 

\begin{lem}
Let $G$ be a definably connected solvable group definable in an o-minimal structure, 
and $A$ a $G$-minimal subgroup of $G$. Then $A\leq Z\o(F(G))$, and 
$C_G(a)=C_G(A)$ for every nontrivial element $a$ in $A$. 
\end{lem}
\proof
By Corollary \ref{GConSolvFGInf}, $A$ has an infinite characteristic abelian definable subgroup. 
Therefore the $G$-minimality of $A$ forces $A$ to be abelian. 
In particular, $A\leq F(G)$. Since $A$ is normal in $F(G)$, Lemma \ref{FactNilpGpsZGInfStrong} 
and the $G$-minimality of $A$ now force that $A\leq Z(F(G))$. Since $A$ is definably connected, 
we have indeed $A\leq Z\o(F(G))$. 

Now $F(G)\leq C_G(A)$, and $G/C_G(A)$ is definably isomorphic to a quotient of $G/F(G)$. 
In particular $G/C_G(A)$ is abelian by Proposition \ref{GConResGFGDivAb}. 
If $A\leq Z(G)$, then clearly $C_G(a)=C_G(A)$ ($=G$) for every $a$ in $A$, and thus we may assume 
$G/C_G(A)$ infinite. Consider the semidirect product $A\rtimes (G/C_G(A))$. 
Since $A$ is $G$-minimal, $A$ is also $G/C_G(A)$-minimal. 
Now an o-minimal version of Zilber's Field Interpretation Theorem 
for groups of finite Morley rank \cite[Theorem 2.6]{MR1779482} applies directly to 
$A\rtimes (G/C_G(A))$. It says that there is an infinite 
interpretable field $K$, with $A\simeq K_+$ and $G/C_G(A)$ an infinite subgroup of $K^{\times}$, and 
such that the action of $G/C_G(A)$ on $A$ corresponds to scalar multiplication. In particular, 
$G/C_G(A)$ acts {\em freely} (or {\em semiregularly} in another commonly used terminology) 
on $A\setminus \{1\}$. This means exactly 
that for any nontrivial element $a$ in $A$, $C_G(a)\leq C_G(A)$, i.e., $C_G(a)=C_G(A)$. 
\qed

\bigskip
The following lemma will be used in a reduction to the solvable case in our proof of 
the main theorem. Of course, it becomes trivial once we now that $G'$ is definable and definably 
connected in definably connected solvable groups (Corollary \ref{CorGoSolvGnGndef} below). 

\begin{lem}\label{LemGResG'InHProper}
Let $G$ be a definably connected solvable group definable in an o-minimal structure. 
If $G$ is not abelian, then $G'\leq H<G$ for some definably connected and definably characteristic 
proper subgroup $H$ definable in $G$. 
\end{lem}
\proof
If $F\o(G)<G$, then it suffices to take $H=F\o(G)$ by Fact \ref{faitFGNilpDef} and 
Proposition \ref{GConResGFGDivAb}. 

Suppose now $G=F\o(G)$, i.e., $G$ nilpotent, and 
assume $G$ nonabelian. Then $G$ is nilpotent of class $n$ for some $n\geq 2$. 
But now by Fact \ref{FactZnGnNilpGps} $G=Z_n(G)$ and $Z_{n-1}(G)<G$, and we also have 
$G'\leq Z_{n-1}(G)$. 
Clearly, all groups $Z_i(G)$ are definable. 
Now one sees as in the proof of Proposition \ref{GConResGFGDivAb} that $G/Z_{n-1}\o(G)$ has a finite 
derived subgroup, which must then be trivial by Fact \ref{ActionConOnFinite} ($b$). 
Hence $G'\leq Z_{n-1}\o(G)<G$, and we may take $H=Z_{n-1}\o(G)$. 
\qed

\begin{cor}\label{CorABSolvAB*}
Let $G$ be a group definable in an o-minimal structure, with $G\o$ solvable. Then $G$ satisfies the 
assumption $(*)$ of Definition \ref{assumption}. 
\end{cor}
\proof
By Lemma \ref{LemGResG'InHProper}, $G\o$ cannot have nonabelian definably simple definable sections. 
\qed

\section{Main theorem}\label{Section8}

We have now enough background to pass to the proof of our main result. 
We first prove a version of Theorem \ref{maintheo} for $A$ and $B$ definably connected. 

\begin{theo}\label{MainTheoConnected}
Let $G$ be a group definable in an o-minimal structure, 
$A$ and $B$ be two definably connected definable subgroups which normalize each other  and 
such that $AB$ satisfies the assumption $(*)$. 
Then $[A,B]$ is a definably connected definable subgroup, and moreover   
$[A,B]=[A,B]_{\dim([A,B])}$ whenever $A$ or $B$ is solvable.  
\end{theo}

\proof
Notice that $[A,B]\leq {A\cap B}$ since $A$ and $B$ normalize each other . 
We proceed by induction on 
$$d:=\min(\dim(A),\dim(B)).$$
Notice that when $d=0$ we have $A$ or $B$ trivial, by definable connectedness, 
and in particular $[A,B]$ is trivial as well. So the induction starts with $d=0$. 
Assume from now on that $G$ is a potential counterexample to our statement, with $d\geq 1$ minimal. 
We start with a series of reductions. Notice that since $AB$ satisfies the assumption $(*)$, 
every definable section of $AB$ does also. 

\begin{claim}\label{ClaimI} 
Let $C$ be a definable subgroup of $A\cap B$ normalized by $A$ and $B$ and with $\dim(C)<d$. 
If $C\nleq{Z(A)\cap Z(B)}$, then $[A,B]$ is definable and definably connected, and 
$[A,B]=[A,B]_{\dim([A,B])}$ whenever $A$ or $B$ is solvable. 
\end{claim}
\proof
Suppose that $C$ does not centralize one of the two groups $A$ or $B$. By symmetry 
we may assume, for instance, $C\nleq Z(B)$. We now have $C\o \nleq Z(B)$ by 
Lemma \ref{FactCentralExtNonCon}. 
In particular $[C\o,B]$ is nontrivial. 
Since $\dim(C)<d$, our inductive assumption implies that $[C\o,B]$ is definable and definably connected, and 
that $[C\o,B]=[C\o,B]_{\dim([C\o,B])}$ whenever $A$ or $B$ is solvable (since $A$ or $B$ solvable implies 
$C\o$ solvable as $C\leq A\cap B$). 

Notice that $[C\o,B]$ is normal in both $A$ and $B$. We now work in $N(U)/U$, where $U=[C\o,B]$, 
and denote by the notation ``$\overline{\phantom{H}}$" the quotients by $U$. 
Since $\dim(U)\geq1$, $\overline{A}$ and 
$\overline{B}$ have dimensions strictly smaller than $\dim(A)$ and $\dim(B)$ respectively. 
Applying our inductive assumption in $N(U)/U$, we get 
$[\overline{A},\overline{B}]$ definable and definably connected, and 
$[\overline{A},\overline{B}]=[\overline{A},\overline{B}]_{\dim([\overline{A},\overline{B}])}$ whenever 
$A$ or $B$ is solvable. 

But clearly 
$$[\overline{A},\overline{B}]=\overline{[A,B]}=[A,B]/U$$ 
since $U\leq [A,B]$, and it follows that $[A,B]$ is definable and definably connected. 

By additivity of the rank, we now get that 
$$\dim([A,B])=\dim(\overline{[A,B]})+\dim(U).$$
Hence, whenever $A$ or $B$ is solvable, we get also that $[A,B]=[A,B]_{\dim([A,B])}$ 
by the additivity of the commutator width provided by Lemma \ref{LemLiftCommutators}. 
\qed

\begin{claim}\label{Claim2}
We may assume ${A\cap B}\nleq {Z(A)\cap Z(B)}$, and $\dim(A\cap B)=d$. 
\end{claim}
\proof
If ${(A\cap B)}\leq{Z(A)\cap Z(B)}$, then in particular $[A,B]\leq{Z(A)\cap Z(B)}$ and 
Corollary \ref{Lem[X,H]definable} would give all the conclusions of Theorem \ref{MainTheoConnected} 
for $[A,B]$. Since we are dealing with a potential counterexample, we may thus assume 
${A\cap B}\nleq {Z(A)\cap Z(B)}$. 

Suppose now $\dim(A\cap B)<d$. Then we may apply Claim \ref{ClaimI} with $C=A\cap B$. 
Since we are now assuming $C\nleq{Z(A)\cap Z(B)}$, it would give all the conclusions of Theorem \ref{MainTheoConnected} for $[A,B]$. As above, we may thus assume $\dim(A\cap B)=d$. 
\qed

\bigskip
By Claim \ref{Claim2}, we now have 
$\dim(A\cap B)=\min(\dim(A),\dim(B))$
and since $A$ and $B$ are definably connected it follows that $A\cap B$ equals $A$ or $B$. By symmetry, 
we may assume $A\cap B=A$, or in other words 
$$A\leq B.$$
Notice that the first part of Claim \ref{Claim2} now says that $A\nleq{Z(A)\cap Z(B)}$. 

\begin{claim}
\label{ClaimWMAA=B}
We may assume $A=B$. 
\end{claim}
\proof
We have to show that if the conclusions of Theorem \ref{MainTheoConnected} are valid for 
$[A,A]=A'$, then they are also valid for $[A,B]$. 
So suppose $[A,A]$ definable and definably connected, and $[A,A]=[A,A]_{\dim(A')}$ 
whenever $A$ is solvable. 

Working modulo $A'$ (which is normal in $B$), we may assume $A$ abelian. 
Then $[A,B]\leq A=Z(A)$ and Corollary \ref{Lem[X,H]definable} gives that 
$[A,B]=[A,B]_{\dim([A,B])}$ is definable and definably connected. 

We can now come back to the original groups $A$ and $B$, i.e., not divided by $A'$. Since 
$A'\leq [A,B]$, we deduce as in the proof of Claim \ref{ClaimI} the definability and definable connectedness 
of $[A,B]$. 
Since $A\leq B$, $A$ is solvable whenever $A$ or $B$ is. Hence we also deduce 
as in the proof of Claim \ref{ClaimI}, using the additivity of the dimension and of the commutator 
width given by Lemma \ref{LemLiftCommutators}, that $[A,B]=[A,B]_{\dim([A,B])}$ whenever 
$A$ or $B$ is solvable. 
\qed

\bigskip
After the preceding series of reductions we are now in the situation 
$$A=B$$
and Claim   \ref{ClaimI},  which of course is still valid, now takes the following form. 
Any proper definable normal subgroup $C$ of $A$ is central in $A$, as otherwise we have finished 
the proof of Theorem \ref{MainTheoConnected}. We split our final analysis in two cases. 

\medskip
\noindent
{\bf Case $A$ solvable. }
Suppose first $A$ solvable. If $A$ is abelian we have of course nothing to do, hence we may suppose 
$A$ nonabelian. By Lemma \ref{LemGResG'InHProper}, $A'\leq C$ for some proper normal 
definably connected definable subgroup $C$, and since we may assume $C\leq Z(A)$ we have 
$A'\leq Z(A)$ (i.e., $A$ is $2$-nilpotent). Now 
Corollary \ref{Lem[X,H]definable} (applied with $H=X=A$) 
gives all the desired conclusions again in this case, including 
the control of the commutator width: $[A,A]=[A,A]_{\dim([A,A])}$. 

\medskip
\noindent
{\bf Case $A$ nonsolvable. } 
Suppose now $A$ non-solvable, i.e., $R(A)<A$. Notice that in this case we just want to conclude 
to the definability and the definable connectedness of $A'$, and we don't care about its commutator width. 

Notice that we actually may suppose with Claim \ref{ClaimI} that $A$ is indecomposable 
in the following sense: it cannot be the product $R_1R_2$ of two proper definable normal subgroups 
$R_1$ and $R_2$. Otherwise, since we may assume $R_1$, $R_2\leq Z(A)$, we would get $A=R_1R_2$ 
abelian, which is excluded at this stage. 

By the structure of semisimple groups definable in o-minimal structures, Fact \ref{FactStrucSemisimpleGps}, 
we get that $A/R(A)$ consists of a single definably simple (and nonabelian) factor. 

Notice also that since we may assume the proper normal definable solvable radical $R(A)$ central in $A$, 
we have also $R(A)\leq Z(A)$, so that $R(A)=Z(A)$. 

To summarize, $A$ is a central extension of a nonabelian definably simple group $A/R(A)$. 

By Fact \ref{FactGDefSimpleGSimple}, $A/Z(A)=[A/Z(A),A/Z(A)]_k$ for some $k$, and lifting commutators 
one gets 
$$A=[A,A]_k\cdot Z(A)$$
and in particular $A=A'\cdot Z(A)$. We distinguish two final subcases. 

\smallskip
\noindent
{\bf Subcase $R(A)=Z(A)$ finite. }
In this case we clearly have ${[A,A]\cap Z(A)}={[A,A]_s\cap Z(A)}$ for some finite $s\geq 1$. 
Now $A'=[A,A]=[A,A]_{k+s}$ by Lemma \ref{LemLiftCommutators}, and in particular $A'$ is definable. 
Since $A=A'Z(A)$ and $Z(A)$ is finite, $A'$ has finite index in $A$. Now since $A$ is definably connected we conclude that $[A,A]=A'=A$ is definable and definably connected, as desired. 

\smallskip
\noindent
{\bf Subcase $R(A)=Z(A)$ infinite. } 
In this final case $A$ is a strict central extension 
of a nonabelian definably simple group. Since our original group $AB$ satisfied the 
assumption $(*)$, the same is true for $A$ at this stage (recall that in our preparatory 
reductions we only took definable sections). Hence, the assumption $(*)$ now says that $A'$ is 
a definable subgroup of $A$. Notice that $A'$ cannot be finite, as otherwise $A$ is abelian 
by Fact \ref{ActionConOnFinite} ($b$). 
If $A'<A$, then our general application of Claim \ref{ClaimI} implies that we may assume 
$A'\leq Z(A)$, which is excluded since $A$ is not solvable. Remains only the case in which 
$[A,A]=A'=A$ is definable and definably connected, as desired. 

\medskip
This finishes the proof of Theorem \ref{MainTheoConnected} in all cases. 
\qed

\bigskip
We now derive from Theorem \ref{MainTheoConnected} a version of 
Theorem \ref{maintheo} with one of the two groups $A$ or $B$ definably connected. 

\begin{cor}\label{CorAoBDefCon}
Let $G$ be a group definable in an o-minimal structure, 
$A$ and $B$ be two definable subgroups which normalize each other  with $A=A\o$ and 
such that $AB\o$ satisfies the assumption $(*)$. 
Then $[A,B]$ is a definably connected definable subgroup, and moreover  
$[A,B]=[A,B]_{\dim([A,B])}$ whenever $A$ or $B\o$ is solvable.  
\end{cor}
\proof
We have $[A,B]$ and $[A,B\o]$ normal in both $A$ and $B$. 

By Theorem \ref{MainTheoConnected} we have $[A,B\o]$ definable and definably connected, 
and $[A,B\o]=[A,B\o]_{\dim([A,B\o])}$ whenever $A$ or $B\o$ is solvable. 

Working in $N([A,B\o])/[A,B\o]$, and denoting by the notation ``$\overline{\phantom{H}}$" the quotients 
by $[A,B\o]$, we may assume that $[\overline{A},\overline{B\o}]=1$. Now the definably connected 
group $\overline A$ acts by conjugation on the finite quotient $\overline{B}/\overline{B\o}$, and 
by Fact \ref{ActionConOnFinite} ($a$) this action must be trivial. This shows that 
$[\overline{A},\overline{B}]\leq \overline{B\o}$. Hence 
$[\overline{A},\overline{B}]\leq {\overline{A}\cap \overline{B\o}}$, and since 
$\overline{A}$ and $\overline{B\o}$ commute we get that 
$$[\overline{A},\overline{B}] \leq Z(\overline{A}).$$
We may now apply Corollary \ref{Lem[X,H]definable} (with $H=\overline{A}$ and $X=\overline{B}$), 
and it gives that $[\overline{A},\overline{B}]$ is a definably connected definable (abelian) 
subgroup of $Z(\overline{A})$, and furthermore that 
$$[\overline{A},\overline{B}]=[\overline{A},\overline{B}]_{\dim([\overline{A},\overline{B}])}.$$
Since $[A,B\o]\leq [A,B]$ and $[\overline{A},\overline{B}]=[A,B]/[A,B\o]$ by lifting of commutators, 
the first point implies that $[A,B]$ is definable and definably connected. 

Now $\dim([A,B])=\dim([\overline{A},\overline{B}])+\dim([A,B\o])$ by additivity of the dimension. 
If $A$ or $B\o$ is solvable, then as in Claim \ref{ClaimI} we get that 
$$[A,B]=[A,B]_{\dim([\overline{A},\overline{B}])+\dim([A,B\o])}=[A,B]_{\dim([A,B])}$$ 
with Lemma \ref{LemLiftCommutators}. 
\qed

\begin{rem}
In Corollary \ref{CorAoBDefCon} we see from the proof that 
$[A,B]/[A,B\o]$ is abelian. In general it need not be trivial, or in other words it is not necessarily 
the case that $[A,B]=[A,B\o]$. Consider for example the following semidirect product 
$B=A\rtimes{\<i\>}$ where $A$ is a divisible abelian group and $i$ is an element of order $2$ 
acting      by inversion on $A$. One may of course realize a group definable in 
an o-minimal structure of this isomorphism type. Clearly $A=B\o$ is definably connected in such a group. 
We have $[A,B\o]=1$, but $[A,B]=A$. 
\end{rem}

We now prove Theorem \ref{maintheo} without assuming any of the two groups $A$ or $B$ 
definably connected. 
 
\bigskip
\noindent
{\bf Proof of Theorem \ref{maintheo}.}
Since $AB$ is a definable subgroup of $G$, we may assume $G=AB$, with $A$ and $B$ 
two normal subgroups. All subgroups $[A,B]$, $[A,B\o]$, and $[A\o,B]$ are normal, in $A$ and $B$, and by 
Corollary \ref{CorAoBDefCon} the two latter subgroups are definable and definably connected. 
By Fact \ref{ProductConGps}, the normal definable subgroup 
$$C:=[A\o,B][A,B\o]$$ 
of $[A,B]$ is definably connected. 

By Corollary \ref{CorAoBDefCon} we also have, as far as the commutator width is concerned, 
$[A\o,B]=[A\o,B]_{\dim([A\o,B])}$ whenever $A\o$ or 
$B\o$ is solvable. Working modulo $[A\o,B]$, we also find similarly as in Claim \ref{ClaimI} or 
Corollary \ref{CorAoBDefCon} (essentially by the additivity of the dimension, 
Lemma \ref{LemLiftCommutators}, and using Corollary \ref{CorAoBDefCon} modulo $[A\o,B]$) 
that the commutator width of $C$ is bounded by $\dim(C)$ 
whenever $A\o$ or $B\o$ is solvable. 

Working modulo $C$, we have $A\o \leq C_A(B)$ and 
$B\o \leq C_B(A)$, and Lemma \ref{CorBaer} gives the finiteness of $[A,B]$ modulo $C$. 
In other words, $[A,B]/C$ is finite, and since $[A,B]$ is a finite extension of the definably connected 
definable subgroup $C$, we have $[A,B]$ definable, and $[A,B]\o=C$. 

This completes our proof of Theorem \ref{maintheo}. 
\qed

\begin{rem}\label{RemComWidthMainTheo}
Clearly, by the last paragraph of the proof of Theorem \ref{maintheo} and Lemma \ref{LemLiftCommutators}, 
we also have $[A,B]$ of finite commutator width in Theorem \ref{maintheo} whenever $A\o$ or $B\o$ 
is solvable. We do not provide explicit bounds 
here, since it depends on the commutator width of the finite section $[A,B]/[A,B]\o$. 
\end{rem}

\bigskip
Applying inductively Theorem \ref{maintheo}, we get the following corollary. 
Recall that it applies to groups $G$ with $G\o$ solvable, by Corollary \ref{CorABSolvAB*}. 

\begin{cor}\label{CorGoSolvGnGndef}
Let $G$ be a group definable in an o-minimal structure such that $G\o$ satisfies the assumption $(*)$. 
Then, for each $n\geq 1$, $G^n$ and $G^{(n)}$ are definable, and definably connected whenever $G$ is. 
\end{cor}

\bigskip
In principle, the commutator widths of the groups $[G^n]\o$ and $[G^{(n)}]\o$ in 
Corollary \ref{CorGoSolvGnGndef} are also controlled by their dimension whenever $G\o$ is solvable. 

We finish this section with two types of questions. 

\begin{question}[The normalization assumption]
\label{QuestionNormAssumption}
Can one prove versions of Theorem \ref{MainTheoConnected}, Corollary \ref{CorAoBDefCon}, 
and Theorem \ref{maintheo}, with $A$ normalizing $B$ only, instead of $A$ and $B$ 
normalizing each other?
\end{question}

A reasonable approach to Question \ref{QuestionNormAssumption} would be to proceed as in the 
proof of Theorem \ref{MainTheoConnected}, but now by induction on $\dim(B)$. Unfortunately, we haven't 
been able to make it work. 

Other natural questions arising concern the control of the commutator width. The first one 
concerns the commutator width in the definably connected case of 
Theorem \ref{MainTheoConnected} (and similarly Corollary \ref{CorAoBDefCon}). 

\begin{question}[The commutator width, I]
\label{QuestionComWidthI}
Can one get rid of the solvability assumption on $A$ or $B$ for the control of the commutator width 
of $[A,B]$ in Theorem \ref{MainTheoConnected}. 
\end{question}

We see from the proof of Theorem \ref{MainTheoConnected} that a positive 
answer to Question \ref{QuestionComWidthI} would only require a uniform control 
on the commutator width of $A'$ in the final case ``$A$ nonsolvable" of that proof. More precisely, 
the question boils down to the following situation. If $A$ is a definably connected group definable 
in an o-minimal structure, with $R(A)=Z(A)$ and $A/R(A)$ (nonabelian) definably simple, and such that 
$A'$ is definable, then can one bound the commutator width of $[A,A]$? This would a priori require a control 
of the commutator width of the definably simple group $A/R(A)$, which is also perfect and of finite commutator 
width by Fact \ref{FactGDefSimpleGSimple}. Here, the problem relates to 
Question \ref{QuestionComWidthDefSimplGp} and has probably a positive answer, possibly with 
a very small commutator width, as commented after Question \ref{QuestionComWidthDefSimplGp}. 
But even if this turned out to be true, one might then 
need to apply Lemma \ref{LemLiftCommutators} to relate the commutator width of $A'$ to that of 
$(A/Z(A))'$, and thus to bound the number $s$ such that $A'\cap Z(A)={[A,A]_s\cap Z(A)}$. 
Even when $Z(A)$ is finite, it is not clear whether the are such uniform bounds (as $A$ varies in the class 
of such central extensions). 

The second natural question about the commutator width concerns the 
``finite bits" modulo the connected components. It can be formulated as follows. 

\begin{question}[The commutator width, II]
\label{QuestionComWidthII}
Fix $G$ a group as in Theorem \ref{maintheo}. Is there a uniform bound on the commutator width 
of the full group $[A,B]$, ``modulo" the commutator width of $[A,B]\o$, as $A$ and $B$ vary as in that 
theorem? 
\end{question}

A positive answer to Question \ref{QuestionComWidthII} would only require to know that, 
for $G$ fixed, the commutator width of the finite section $[A,B]/[A,B]\o$ is uniformly 
bounded as long as $A$ and $B$ vary in the set of 
all definable subgroups of $G$ satisfying the assumption $(*)$ 
(Remark \ref{RemComWidthMainTheo}). 
By the specific structure of $G$ following Fact \ref{FactStrucSemisimpleGps} 
($G\o/R(G\o)$ is a direct product of definably simple groups) 
and a closer look at finite sections of the solvable group $R(G\o)$, it might be true. 

\section{Further questions and remarks}\label{Section9}

In any group $G$ of finite Morley rank we have the following: if $H$ is any definable and 
definably connected subgroup, and $X$ an arbitrary subset of $G$, then 
$[X,H]$ is definable and definably connected \cite[Corollary 5.29]{MR1321141}. 
Of course the proof of this general statement 
uses the definition of the Morley rank in a crucial way via Zilber's stabilizer argument, 
but one may wonder when such stronger versions of Corollary \ref{CorAoBDefCon} 
remain valid in the context of groups definable in o-minimal structures. 

\begin{question}\label{question[H,X]defcon}
In a group definable in an o-minimal structure, let $H$ be a definable and 
definably connected subgroup and $X$ an arbitrary 
subset. When is it true that $[H,X]$ is definable and definably connected?
\end{question}

Corollary \ref{Lem[X,H]definable} remains our most general answer to 
Question \ref{question[H,X]defcon}. On the other hand, Question \ref{question[H,X]defcon} 
may fail even if $X$ is a normal subgroup of $H$, as the following example shows. 

\begin{example}\label{Example[H,X]NotDef}
Let $H$ be a definably compact definably connected definably simple nonabelian group definable 
in an o-minimal structure. Taking a sufficiently saturated model, we know that the 
smallest type-definable subgroup $X=H^{\circ \circ}$ of bounded index in $H$ is 
normal, proper, and nontrivial (by \cite[2.12]{MR2114965}, and 
\cite[3.6]{MR2114965} or \cite{MR2139910}). Then $[H,X]$ is not definable. 
Otherwise, the proper normal definable subgroup $[H,X]\leq X<H$ would be trivial by definable simplicity 
of $H$, implying the existence of a nontrivial center in $H$. 
\end{example}

We manage however to treat the case $H\subseteq X$ as a 
corollary of Theorem \ref{MainTheoConnected} as follows. 

\begin{cor}
\label{CorFinalHX}
Let $G$ be a group definable in an o-minimal structure, $H$ a definably connected definable 
subgroup which satisfies the assumption $(*)$, and $X$ any subset such that 
$H\subseteq X\subseteq N(H)$. Then $[H,X]$ is a definable and definably connected subgroup of $H$. 
\end{cor}
\proof
The group $H'$ is definable and definably connected by Theorem \ref{MainTheoConnected}. 
Since $H'$ is normalized by $H$ and $X$, we may replace $G$ by 
$N(H')$, and working modulo $H'$ we may assume $H$ abelian. 
Now, since $[H,X]\leq H=Z(H)$, Corollary \ref{Lem[X,H]definable} gives the definability and 
definable connectedness of $[H,X]$. 
\qed

\bigskip
Of course, we can also obtain $[H,X]=[H,X]_{\dim([H,X])}$ as in the proof of Claim \ref{ClaimWMAA=B} 
when $H$ is solvable in Corollary \ref{CorFinalHX}. 

\bibliographystyle{plain}
\bibliography{biblio}

\end{document}